\documentclass[]{aiaa-tc}% insert '[draft]' option to show overfull boxes

\usepackage{aiaa_packages}

 \title{Low-Thrust Trajectory Design Using Reachability Sets near Asteroid 4769 Castalia}

 \author{
  Shankar Kulumani\thanksibid{1}%
    \thanks{Doctoral Student, \href{mailto:skulumani@gwu.edu}{skulumani@gwu.edu}. Student AIAA Member.}
  \ and Taeyoung Lee\thanksibid{2}\thanks{Associate Professor, \href{mailto:tylee@gwu.edu}{tylee@gwu.edu}.}\\
  {\normalsize\itshape
   George Washington University, Washington, District of Columbia }\\
   }

 \AIAApapernumber{YEAR-NUMBER}
 \AIAAconference{Conference Name, Date, and Location}
 \AIAAcopyright{\AIAAcopyrightD{YEAR}}

\begin{document}

\maketitle

\begin{abstract}\label{sec:abstract}

We present a computational approach for the design of continuous low thrust transfers around an asteroid.
These transfers are computed through the use of a reachability set generated on a lower dimensional \Poincare surface.
Complex, long duration transfer trajectories are highly sensitive to the initial guess and generally have a small region of convergence.
Computation of the reachable set alleviates the need to generate an accurate initial guess for optimization.
From the reachable set, we chose a trajectory which minimizes a distance metric towards the desired target.
Successive computation of the reachable set allows for the design of general transfer trajectories which iteratively approach the target.
We demonstrate this method by determining a transfer trajectory about the asteroid 4769 Castalia.

\end{abstract}

\section{Introduction}\label{sec:introduction}

% Motivation for missions/studying asteroids
Small solar system bodies, such as asteroids and comets, are of significant interest to the scientific community; these small bodies offer great insight into the early formation of the solar system.
This insight offers additional detail into the formation of the Earth and also the probable formation of other planetary bodies.
Of particular interest are those near-Earth asteroids (NEA) which inhabit heliocentric orbits in vicinity of the Earth.
These easily accessible bodies provide attractive targets to support space industrialization, mining operations, and scientific missions.
NEAs potentially contain many materials such as those useful for propulsion, construction, or for use in semiconductors.
Also, many bodies contain highly profitable materials, such as precious or strategic metals~\cite{ross2001}.
In addition, these NEAs are also of concern for their potential to impact the Earth.
Asteroids and comets are the greatest threat to future civilizations and as a result there is a focused effort to mitigate these risks~\cite{wie2008}.
In spite of the great interest in asteroids, the operation of spacecraft in their vicinity is a challenging problem.

% Difficulty in system model
While there has been significant study of interplanetary transfer trajectories, relatively less analysis has been conducted on operations in the vicinity of asteroids.
The dynamic environment around asteroids is strongly perturbed and challenging for analysis and mission operation~\cite{scheeres1994,scheeres2000}.
Due to their low mass, which results in a low gravitational attraction, asteroids may have irregular shapes and potentially chaotic spin states.
Furthermore, since the magnitude of the gravitational attraction is relatively small, non-gravitational effects, such as solar radiation pressure or third-body effects, become much more significant.
As a result, the orbital environment is generally quite complex and it is difficult to generate analytical insights.

% most use a spherical harmonic model or a ellipsoid model but we use a polyhedron model
An accurate gravitational potential model is necessary for the operation of spacecraft about asteroids.
Additionally, a detailed shape model of the asteroid is needed for trajectories passing close to the body.
The classic approach is to expand the gravitational potential into a harmonic series and compute the series coefficients.
Radio tracking data of an orbiting spacecraft allows one to estimate the series coefficients.
The harmonic representation is guaranteed to converge outside of the circumscribing sphere and can be truncated at a finite order based on accuracy requirements~\cite{scheeres2012a}.
However, the harmonic expansion is always an approximation as a result of the infinite order series used in the representation.
Additionally, the harmonic model used outside of the circumscribing sphere is not guaranteed to converge inside the sphere.

A popular approach to deal with this divergence is to use a different gravitational model within the circumscribing sphere.
For example, Reference~\citenum{scheeres2000} uses both a polyhedron field and a spherical harmonic expansion to represent the full gravitational field of the body.
The spherical harmonic coefficients are computed from the polyhedron shape model using a constant density assumption~\cite{werner1997}.
This model is applied close to the body while the harmonic expansion is applied when outside the circumscribing sphere.
Similarly, Reference~\citenum{herrera2014} uses two different harmonic expansion models to represent the gravitational field inside and outside of the circumscribing sphere.
In this case, the coefficients must be carefully chosen to ensure continuity of the gravitational field at the circumscribing sphere.
In both cases, this type of approach results in a cumbersome gravity field expression that requires additional constraints to ensure continuity and validity at the radius of the circumscribing sphere.
Any integration software will need to incorporate a switching mechanism between the models when crossing the circumscribing sphere.
Instead, we use a different approach and model the asteroid as a constant-density polyhedron.

The polyhedron model results in an exact closed form expression of the gravitational potential field~\cite{werner1994,werner1996}.
This type of model results in the exact potential up to the accuracy of the shape model and a constant density assumption.
However, the calculation of the potential, or the acceleration, requires a summation over every face of the polyhedron. 
As a result, it typically requires a large amount of computation in contrast to the harmonic expansion models. 
However, the formulation is well suited to parallelization and improvements with efficient coding practices. 
Finally, the polyhedron method is well suited to trajectories passing close to the body and offers a simple metric to determine if a particle is inside the asteroid.
The use of the polyhedron model results in a single expression of the gravitational field that is valid everywhere around the body.

% use of low thrust propulsion
The application of optimal control methods for orbital trajectory design is nontrivial.
Frequently, insight into the problem or intuition on the part of the designer is required to determine initial conditions that will converge to the optimal solution.
However, the asteroid system dynamics are nonlinear and exhibit chaotic behaviors.
This makes solving the optimization problem highly dependent on the initial condition.
Similar to the three-body problem, there is an insufficient number of analytical constants to derive an analytical solution in general.
As a result, accurate numerical methods are required to determine optimal solutions.
These methods are critically dependent on an accurate initial guess in order to allow for convergence.

% Discuss previous work and their drawbacks to determine transfers about asteroids
We model the motion of particles around asteroids using the restricted full two-body problem.
The dynamics of a spacecraft about small bodies is very similar to that of the three-body problem.
This model has many similarities with the restricted three body problem, and much of the theory developed for the three-body problem is also applicable~\cite{mondelo2010,herrera2014}.
In addition, there has been a large amount of work on the optimal control of spacecraft orbital transfers in the three-body problem~\cite{mingotti2011,grebow2011}.
Typically, the optimal control problem is solved via direct methods, which approximate the continuous time problem as a parameter optimization problem.
The state and control trajectories are discretly parameterized and solved in the form of a nonlinear optimization problem.
Alternatively, indirect methods apply calculus of variations to derive the necessary conditions for optimality. 
This yields a lower dimensioned problem as compared to the direct approach.
In addition, satisfication of the necessary conditions guarantee local optimality in contrast to direct methods which result in sub-optimal solutions.

In this paper, we extend the design method previously developed in the three-body problem to motion about asteroids~\cite{kulumani2015}.
Our systematic approach avoids the difficulties in selecting an appropriate initial guess for optimization.
We instead utilize the concept of the reachability set to enable a simple methodology of selecting initial conditions to achieve general orbital transfers.
This method allows the spacecraft to depart from fixed periodic solutions through the use of a low-thrust propulsion system.
In addition, we utilize a polyhedron gravitational model which is accurate and is globally applicable about the asteroid.

We formulate an optimal control problem to calculate the reachability set on a lower dimensional \Poincare section.
Given an initial condition and fixed time horizon, the reachable set is the set of states attainable, subject to the operational constraints of the spacecraft.
The generation of the reachable set allows for a more systematic method of determining initial conditions and eases the burden on the designer.
The \Poincare section reduces the dimensionality of the system dynamics to the study of a related discrete update map.
This allows for the design of complex transfer trajectories on a lower dimensional space.
Rather than relying on intuition or insight into the problem, trajectories are chosen which minimize a distance metric toward a desired target on the \Poincare section.
This simple methodology allows for extended transfer trajectories which iteratively approach a desired target orbit.

In short, the authors present a systematic method of generating optimal transfer orbits about asteroids.
Typically, optimal transfers are generated using a direct optimization method which results in a sub-optimal solution.
This paper present an indirect optimal control formulation to generate the reachability set on a \Poincare section.
Using the reachability set on the \Poincare section allows for a simple method of choosing trajectories which approach a target.
In addition, the reachability set gives an indication of the regions of the phase space accessible to the spacecraft.
This method allows us to avoid the difficulties inherent in choosing valid initial conditions for the computation of optimal transfer trajectories.
We develop the optimal control formulation and apply this method to a transfer about asteroid 4769 Castalia.

\section{Asteroid Model}\label{sec:asteroid_model}

In this analysis we consider transfers about the asteroid 4769 Castalia.
Castalia has an accurate shape model and is also considered a potentially hazardous asteroid with a possibility of Earth impact~\cite{hudson1994}.
We model the gravitational potential field of Castalia using a polyhedron gravity model instead of using a spherical harmonic expansion.
The spherical harmonic expansion is a popular method of representing the gravity field~\cite{scheeres1996}.
Approximations are possible by truncating the infinite order series to fixed set of coefficients, with the most important terms corresponding to the second order and degree~\cite{scheeres1994}. 
However, when evaluated close to the body the series expansion will diverge and is no longer accurate. 
Therefore, the spherical harmonic representation is not ideal for landing trajectories or those passing close to the surface.

% benefits of the polyhedron model
A polyhedral model of the surface of an asteroid can be determined from remote optical or radar sensors.
The faces of the polyhedron can be large or small and allow for fine detail such as depression, craters, ridges, or interior voids.
In addition, there is no requirement for the body to be modeled at a uniformly high resolution so small details can be incorporated with minimal cost.
From the shape model, an analytical, closed form expression for the gravitational potential can be derived.
The polyhedral approach provides an accurate gravitational model consistent with the resolution of the shape and the chosen discretization.
Furthermore, the polyhedron model is an exact solution up to the surface of the body.
Therefore, this model is ideal for missions traversing large regions both close and far from the asteroid.

\subsection{Polyhedron Gravity Model}\label{sec:polyhedron_model}

We represent the gravitational potential of the asteroid using a polyhedron gravitation model.
This model is composed of a polyhedron, which is a three-dimensional solid body, that is defined by a series of vectors in the body-fixed frame.
The vectors define vertices in the body-fixed frame as well as planar faces which compose the surface of the asteroid.
We assume that each face is a triangle composed of three vertices and three edges.
As a result, only two faces meet at each edge while three faces meet at each vertex.
Only the body-fixed vectors, and their associated topology, is required to define the exterior gravitational model.
References~\citenum{werner1994} and~\citenum{werner1996} give a detailed derivation of the polyhedron model.
Here, we summarize the key developments and equations required for implementation.

Consider three vectors \( \vecbf{v}_1, \vecbf{v}_2, \vecbf{v}_3 \in \R^{3 \times 1} \), assumed to be ordered in a counterclockwise direction, which define a face.
It is easy to define the three edges of each face as
\begin{align}\label{eq:edges}
    \vecbf{e}_{i+1,i} = \vecbf{v}_{i+1} - \vecbf{v}_i \in \R^{3 \times 1 },
\end{align}
where the index \( i \in \parenth{1,2,3} \) is used to permute all edges of each face.
Since each edge is a member of two faces, there exist two edges which are defined in opposite directions between the same vertices.
We can also define the outward normal vector to face \( f\)  as
\begin{align}\label{eq:face_normal}
    \hat{\vecbf{n}}_f &= \parenth{\vecbf{v}_{2} - \vecbf{v}_1} \times \parenth{\vecbf{v}_{3} - \vecbf{v}_2} \in \R^{3 \times 1},
\end{align}
and the outward facing normal vector to each edge as
\begin{align}\label{eq:edge_normal}
    \hat{\vecbf{n}}_{i+1,i}^f &= \parenth{\vecbf{v}_{i+1} - \vecbf{v}_i} \times \hat{\vecbf{n}}_f \in \R^{3 \times 1}.
\end{align}
For each face we define the face dyad \( \vecbf{F}_f \) as
\begin{align}\label{eq:face_dyad}
    \vecbf{F}_f &= \hat{\vecbf{n}}_f \hat{\vecbf{n}}_f \in \R^{3 \times 3}.
\end{align}
Each edge is a member of two faces and has an outward pointing edge normal vector, given in~\cref{eq:edge_normal}, perpendicular to both the edge and the face normal.
For the edge connecting the vectors \( \vecbf{v}_1 \) and \( \vecbf{v}_2 \), which are shared between the faces \(A\) and \( B\), the per edge dyad is given by
\begin{align}\label{eq:edge_dyad}
    \vecbf{E}_{12} = \hat{\vecbf{n}}_A \hat{\vecbf{n}}_{12}^A + \hat{\vecbf{n}}_B \hat{\vecbf{n}}_{21}^B \in \R^{3 \times 3}.
\end{align}
The edge dyad \( \vecbf{E}_e  \), is defined for each edge and is a function of the two adjacent faces meeting at that edge.
The face dyad \( \vecbf{F}_f \), is defined for each face and is a function of the face normal vectors.

Let \( \vecbf{r}_i \in \R^{3 \times 1} \) be the vector from the spacecraft to the vertex \( \vecbf{v}_i \) and it's length is given by \( r_i = \norm{\vecbf{r}_i} \in \R^{1} \).
The per-edge factor \( L_e \in \R^{1}\), for the edge connecting vertices \( \vecbf{v}_i \) and \( \vecbf{v}_j \), with a constant length \( e_{ij} = \norm{\vecbf{e}_{ij}} \in \R^1\) is
\begin{align}\label{eq:edge_factor}
    L_e &= \ln \frac{r_i + r_j + e_{ij}}{r_i + r_j - e_{ij}}.
\end{align}
For the face defined by the vertices \( \vecbf{v}_i, \vecbf{v}_j, \vecbf{v}_k \) the per-face factor \( \omega_f \in \R^{1} \) is
\begin{align}\label{eq:face_factor}
    \omega_f &= 2 \arctan \frac{\vecbf{r}_i \cdot \vecbf{r}_j \times \vecbf{r}_k}{r_i r_j r_k + r_i \parenth{\vecbf{r}_j \cdot \vecbf{r}_k} + r_j \parenth{\vecbf{r}_k \cdot \vecbf{r}_i} + r_k \parenth{\vecbf{r}_i \cdot \vecbf{r}_j}} .
\end{align}
The gravitational potential due to a constant density polyhedron is given as
\begin{align}\label{eq:potential}
    U(\vecbf{r}) &= \frac{1}{2} G \sigma \sum_{e \in \text{edges}} \vecbf{r}_e \cdot \vecbf{E}_e \cdot \vecbf{r}_e \cdot L_e - \frac{1}{2}G \sigma \sum_{f \in \text{faces}} \vecbf{r}_f \cdot \vecbf{F}_f \cdot \vecbf{r}_f \cdot \omega_f \in \R^1,
\end{align}
where \( \vecbf{r}_e\) and \(\vecbf{r}_f \) are the vectors from the spacecraft to any point on the respective edge or face, \( G\) is the universal gravitational constant, and \( \sigma \) is the constant density of the asteroid.
Furthermore we can use these definitions to define the attraction, gravity gradient matrix, and Laplacian as
\begin{align}
    \nabla U ( \vecbf{r} ) &= -G \sigma \sum_{e \in \text{edges}} \vecbf{E}_e \cdot \vecbf{r}_e \cdot L_e + G \sigma \sum_{f \in \text{faces}} \vecbf{F}_f \cdot \vecbf{r}_f \cdot \omega_f \in \R^{3 \times 1} , \label{eq:attraction}\\
    \nabla \nabla U ( \vecbf{r} ) &= G \sigma \sum_{e \in \text{edges}} \vecbf{E}_e  \cdot L_e - G \sigma \sum_{f \in \text{faces}} \vecbf{F}_f \cdot \omega_f \in \R^{3 \times 3}, \label{eq:gradient_matrix}\\
    \nabla^2 U &= -G \sigma \sum_{f \in \text{faces}}  \omega_f \in \R^1 .\label{eq:laplacian}
\end{align}

One interesting thing to note is that both~\cref{eq:face_dyad,eq:edge_dyad} can be precomputed without knowledge of the position of the satellite.
They are both solely functions of the vertices and edges of the polyhedral shape model and are computed once and stored.
Once a position vector \( \vecbf{r} \) is defined, the scalars given in~\cref{eq:edge_factor,eq:face_factor} can be computed for each face and edge.
Finally,~\cref{eq:potential} is used to compute the gravitational potential on the spacecraft.
The Laplacian, defined in~\cref{eq:laplacian}, gives a simple method to determine if the spacecraft has collided with the body~\cite{werner1996}. 

\begin{figure}
    \centering
    \begin{subfigure}[htbp]{0.45\textwidth}
        \includegraphics[width=\textwidth]{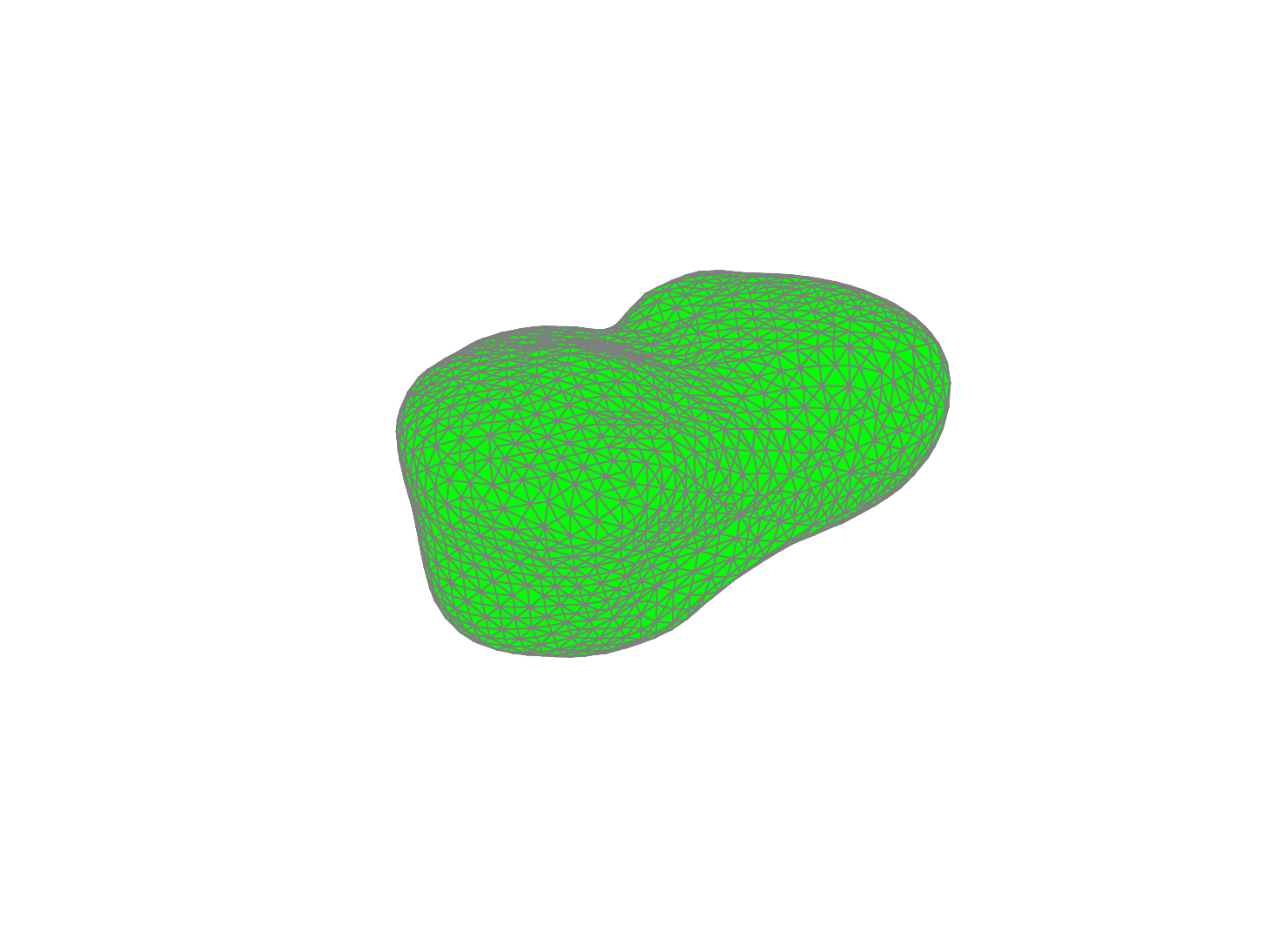}
        \caption{3D Shape Model of Castalia} \label{fig:castalia_3d}
    \end{subfigure}~ %add desired spacing between images, e. g. ~, \quad, \qquad, \hfill etc. %(or a blank line to force the subfigure onto a new line)
    \begin{subfigure}[htbp]{0.45\textwidth}
        \includegraphics[width=\textwidth]{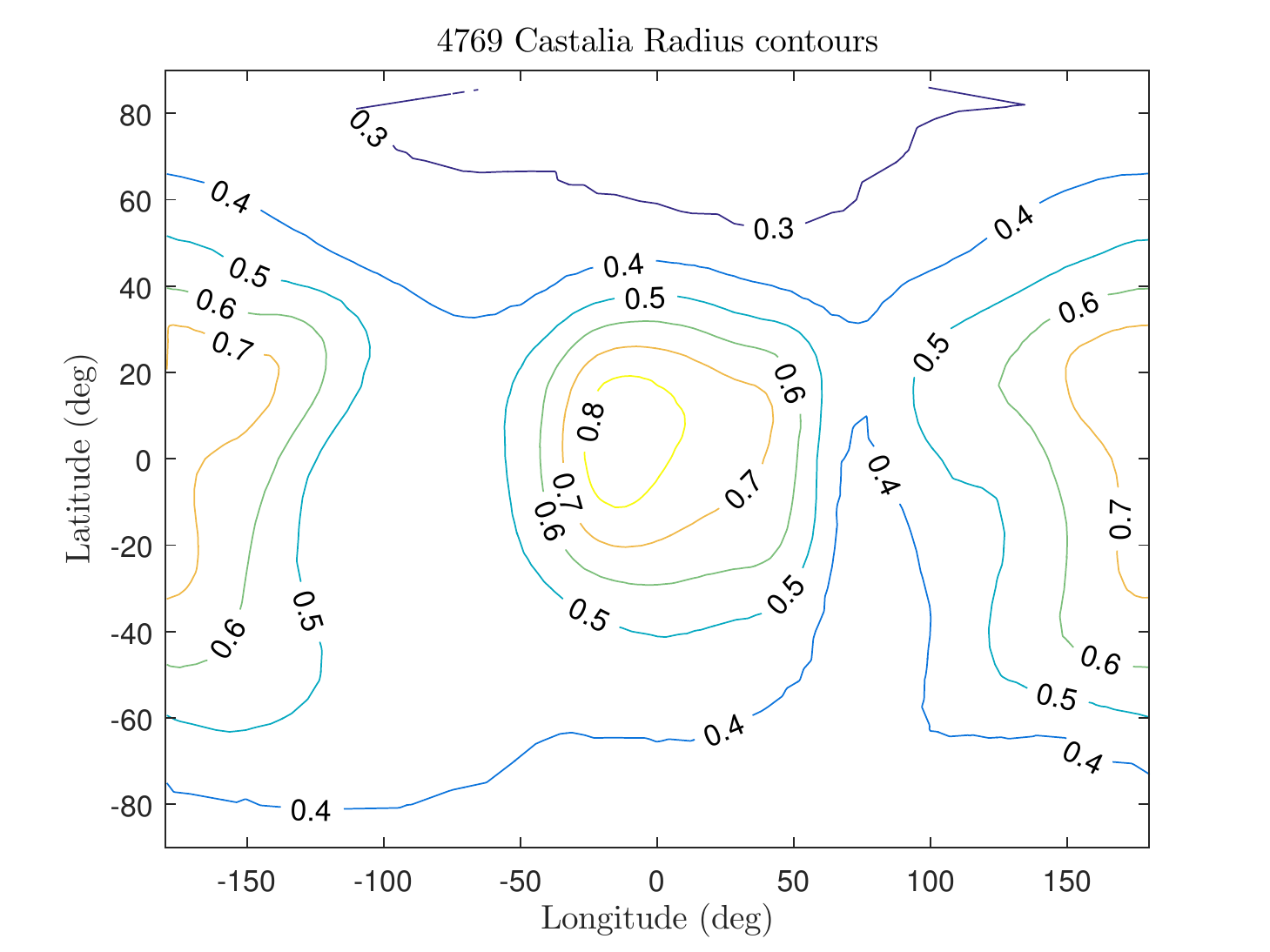}
        \caption{Radius contours of Castalia} \label{fig:radius_contour}
    \end{subfigure} ~ %add desired spacing between images, e. g. ~, \quad, \qquad, \hfill etc. %(or a blank line to force the subfigure onto a new line)
    \caption{Polyhedron Shape Model of 4769 Castalia}
    \label{fig:castalia}
\end{figure}

In this work, we consider trajectories about asteroid 4769 Castalia.
Doppler radar images, obtained at the Arecibo Observatory in 1989, are used to determine a shape model of Castalia~\cite{hudson1994,neese2004}.
We use the estimated rotation period of \SI{4.07}{\hour} with a nominal density of \SI{2.1}{\gram\per\centi\meter\cubed}~\cite{scheeres1996}.
The shape model is composed of \num{4092} triangular faces and a rendering of the asteroid is provided in~\cref{fig:castalia_3d}.
In addition, we show a contour plot of the radius of Castlia in~\cref{fig:radius_contour}.

\subsection{Spacecraft Equations of Motion}\label{sec:sc_eoms}

The motion of a massless particle, or spacecraft, about an asteroid shares many similarities with that of the three-body problem.
As is typical in the three-body problem, the equations of motion are usually represented in a uniformly rotating frame aligned with the two primaries.
Similarly, the equations of motion about an asteroid are also defined in a body-fixed frame with uniform rotation.
In this reference frame, the gravitational potential field is time invariant and only a function of the position of the particle.
In addition, since the rotational rate of the asteroid is constant, the equations of motion are time invariant.
Finally, the use of the rotating reference frame allows for much greater insight into the dynamic structure of the behavior around the asteroid.

We define a reference frame originating at the center of mass of the asteroid.
The body-fixed reference frame is composed of the unit vectors \( \hat{\vecbf{x}} , \hat{\vecbf{y}}, \hat{\vecbf{z}} \), which are aligned along the principal axes of smallest, intermediate, and largest moment of inertia, respectively.
The body-fixed equations of motion of a massless particle about an arbitrarily rotating asteroid are given by
\begin{align}\label{eq:body_eoms}
    \ddot{\vecbf{r}} + 2 \vecbf{\Omega} \times \dot{\vecbf{r}} + \vecbf{\Omega} \times \parenth{ \vecbf{\Omega} \times \vecbf{r} } + \dot{\vecbf{\Omega}} \times \vecbf{r} = \nabla U(\vecbf{r}) ,
\end{align}
where \( \vecbf{\Omega} \in \R^{3 \times 1}\) is the instantaneous angular velocity vector of the asteroid represented in the body-fixed frame, \( \vecbf{r} \) is the position of the particle in the body-fixed frame, and \( \nabla U(\vecbf{r}) \) is the gradient of the gravitational potential~\cite{scheeres2012a}.
We assume that the asteroid rotates at a uniform rate, \( \norm{\vecbf{\Omega}} = \omega \in \R^1 \), about the axis of the maximum moment of inertia, i.e.\ \( \vecbf{\Omega} = \omega \hat{\vecbf{z} }\).
As a result, we can represent the equations of motion in scalar form as
\begin{align} \label{eq:eoms}
    \begin{split}
        \ddot{x} - 2 \omega \dot{y} - \omega^2 y &= U_x , \\
        \ddot{y} + 2 \omega \dot{x} - \omega^2 x &= U_y , \\
        \ddot{z} &= U_z .
    \end{split}
\end{align}

In this situation, the state is defined as \( \vecbf{x} = \bracket{\vecbf{r}~\>\vecbf{v}}^T \in \R^{6 \times 1}\) with \(\vecbf{r} = \bracket{x~\>y~\>z}^T \in \R^{3\times1}\) and \(\vecbf{v}= \bracket{ \dot{x}~\>\dot{y}~\>\dot{z} }^T \in \R^{3\times1}\) representing the position and velocity with respect to the body-fixed frame, respectively.
We further assume that our spacecraft is capable of exerting a translational acceleration, \( \vecbf{u} \in \R^{3\times1} \), in any direction, while subject to a maximum magnitude constraint, \( \norm{\vecbf{u}} \leq u_m \).
This is typical of many spacecraft which offer full rotational freedom and can direct a potentially varying force or acceleration in any direction.
The equations of motion may be rewritten in state space form as
\begin{align}\label{eq:state_space_eoms}
    \begin{bmatrix} \dot{\vecbf{r}} \\ \dot{\vecbf{v}} \end{bmatrix} &=
    \begin{bmatrix}\vecbf{v} \\ \vecbf{g} \parenth{\vecbf{r}} + \vecbf{h}\parenth{\vecbf{v}} + \vecbf{u} \end{bmatrix} ,
\end{align}
where the terms \(\vecbf{g} \parenth{\vecbf{r}} \) and \( \vecbf{h}\parenth{\vecbf{v}} \) are given by
\begin{align}\label{eq:state_space_terms}
    \vecbf{g}\parenth{\vecbf{r}} = \begin{bmatrix}  U_x + \omega^2 x \\ U_y + \omega^2 y \\ U_z \end{bmatrix} ,\quad
    \vecbf{h}\parenth{\vecbf{r}} = \begin{bmatrix} 2 \omega \dot{y} \\ -2 \omega \dot{x} \\ 0 \end{bmatrix} .
\end{align}

Since Castalia is a uniformly rotating asteroid, the equations of motion are time invariant when represented in the body-fixed frame.
In addition, there exists an integral of motion, or a conserved quantity, that is constant for all motion of a particle.
The Jacobi constant, \( J (\vecbf{r} , \vecbf{v} ) \), is given by
\begin{align}\label{eq:jacobi}
    J \parenth{\vecbf{r}, \vecbf{v}} = \frac{1}{2} \omega^2 \parenth{x^2 + y^2} + U(\vecbf{r}) - \frac{1}{2} \parenth{\dot{x}^2 + \dot{y}^2 + \dot{z}^2} .
\end{align}
The Jacobi constant functions in a similar manner as used in three-body problem~\cite{szebehely1967}.
We can define zero-velocity surfaces using the Jacobi constant by fixing the value to a desired constant.
The zero-velocity surfaces are the locus of points where the kinetic energy and hence velocity vanishes.
Just as in the three-body problem, the Jacobi constant in~\cref{eq:jacobi} divides the phase space into distinct realms of possible motion.
Similarly, there exist, in general, four equilibrium points and also their associated stable and unstable manifolds~\cite{scheeres1996,scheeres1994}.
The properties of these manifolds play a critical role in the dynamics of trajectories in their vicinity.

\section{Reachability Set on a \Poincare Section}\label{sec:reachability}

Typical optimal control methods, including both indirect and direct based methods, are highly dependent on an accurate initial guess.
For indirect optimization, which is based on the calculus of variations, this results in the well-known two-point boundary value problem.
Insight into the problem or insight by the designer is usually required to determine appropriate initial costates that will converge to the optimal solution and satisfy the desired constraints.
To avoid this issue, we utilize the concept of the reachability set on a lower dimensional \Poincare section.
By repeatedly constructing the reachability set, we can achieve general transfers by determining set intersections on the \Poincare section.
This alleviates the need to determine an accurate initial guess while offering some insight into the dynamics of neighboring trajectories

The reachable set contains all possible trajectories that are achievable over a fixed time horizon from a defined initial condition, subject to the constraints of the system.
Reachability theory has been applied to collision avoidance and safety planning in aerospace systems~\cite{holzinger2009,holzinger2011b}.
The theory supporting reachability analysis is directly derivable from optimal control theory~\cite{varaiya2000,lygeros2004}.
Analytic computation of reachability sets is only possible for a small class of potential systems.
Here, we use numerical methods to solve a related optimal control problem, which approximates a single solution that lies on the reachable set.
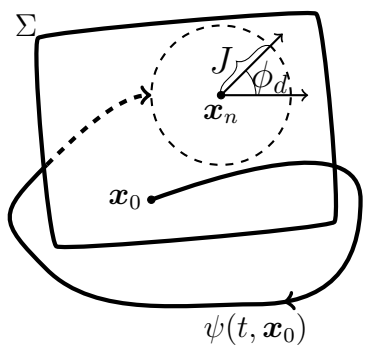
\begin{figure}
    \centering
    \begin{scaletikzpicturetowidth}{0.3\textwidth}
    \begin{tikzpicture}[scale=\tikzscale]
        \coordinate [label=left:\textcolor{black}{\large \(\vecbf{x}_0\)}] (x0) at (-1,-2);
        \coordinate [label=below:\textcolor{black}{\large  \(\vecbf{x}_n\)}] (xn) at (1,1);
        \coordinate [label=left:\textcolor{black}{\large  \(\Sigma\)}] (sigma) at (-4,3);
        %\coordinate [label=below:\textcolor{black}{\large  \(P(\vec{x})\)}] (P) at (0,-3.5);
        % define the path of the flow with coordinates
        \coordinate [label=right:\textcolor{black}{}] (f1) at (5,-2);
        \coordinate [label=below:\textcolor{black}{\large  \(\psi(t,\vecbf{x}_0)\)}] (f2) at (2,-5);
        \coordinate [label=right:\textcolor{black}{}] (f3) at (-4,-4);
        \coordinate [label=right:\textcolor{black}{}] (f4) at (-4,-1);
        
    %   \draw[help lines] (-10,-10) grid (10,10); %grid
        \filldraw [black] (x0) circle [radius=3pt];
        \filldraw [black] (xn) circle [radius=3pt];
    
        \draw [ultra thick,black,->-](x0) to[out=20,in=90,distance=2cm] (f1) to[out=-90,in=0,distance=2cm] (f2) to[out=180,in=-45,distance=2cm] (f3) to[out=135,in=-135,distance=2cm] (f4) ;
        \draw [ultra thick, black,dashed,->] (f4) to[out=45,in=180,distance=1cm] ($(xn)-(2,0)$);
        
        \draw [ultra thick] plot [smooth cycle, tension=0.1, rotate=5] coordinates { (-4,-3) (4,-3) (4,3) (-4,3) };
    
        \draw [thick,dashed] (xn) circle [radius=2cm]; % reachability set
    
        \draw [thick,->] (xn) -- ($(xn) + (2.5,0)$);
        \draw [thick,rotate=45,->] (xn) -- ($(xn) + (2.5,0)$);
        \draw ($(xn) + (1,0)$) arc [start angle=0,end angle=45, radius=1];
        \node [draw=none] at (2.4,1.5) {\Large \(\phi_d\)};
        \draw [decorate,decoration={brace,amplitude=5pt},rotate=45] (xn) -- ($(xn) + (2,0)$);
        \node [draw=none] at ($ (xn) + (0,1) $) {\Large \( J \)};
    \end{tikzpicture}
    \end{scaletikzpicturetowidth}
    \caption{Reachability set on a \Poincare section\label{fig:reachability_set}}
\end{figure}

We seek to approximate the reachability set on a \Poincare section by solving an optimal control problem.
The \Poincare section is chosen in a manner similar to the previous work in both the three-body problem as well analysis performed around asteroids to determine periodic orbits.
Typically, analysis for the three-body problem relies heavily on symmetries in the force fields.
However, in our system model, the gravitational potential, given by~\cref{eq:potential}, has no symmetries.
In spite of this, it is still possible to determine periodic solutions through the application of a \Poincare map with the surface of section chosen normal to a surface in the phase space~\cite{scheeres2000}.
For a periodic orbit, the trajectories will intersect the \Poincare section at two distinct points every half orbit.
With the addition of a low thrust control input, we are able to expand the reachable set from a distinct point to a larger area on the \Poincare section.
\Cref{fig:reachability_set} illustrates this methodology.
Without any control input, the trajectories will follow the system dynamics, \( \psi(t, \vecbf{x}_0 ) \) and intersect the \Poincare section at \( \vecbf{x}_n\).
The addition of a control input allows the spacecraft to depart from the natural dynamics and intersect the section at another location denoted by the dashed circle.
We use the cost function \( J \) to define a distance metric on the \Poincare section.
Maximization of \( J \), or the minimization of \( -J \), along various directions, which are parameterized using \( \phi_d \), on the \Poincare section allows us to generate the largest reachability set under the bounded control input.

We define the \Poincare section as the surface normal to \( y = 0 \).
Following convention, the \Poincare map is defined as the map from one transversal crossing of the surface \( y = 0\) to the next.
Using the method of Reference~\citenum{scheeres2000}, we remove \( y \) and \( \dot{y} \) from consideration and create a four-dimensional map.
The \Poincare section, represented by \( \Sigma \), then becomes
\begin{align}\label{eq:poincare_section}
    \Sigma = \braces{\parenth{x, \dot{x}, z, \dot{z}} | y(t_f) = 0 }.
\end{align}
We use this section to compute periodic orbits that serve as the initial and target states of our transfer.
In addition, this section serves as a lower dimensional space upon which we approximate the reachability set.

An optimal control problem is defined by the cost function
\begin{align}\label{eq:cost}
    J = -\frac{1}{2} \left( \vecbf{x}(t_f) - \vecbf{x}_{n}(t_f)\right)^T 
    Q
    \left( \vecbf{x}(t_f) - \vecbf{x}_{n}(t_f)\right) ,
\end{align}
where \( \vecbf{x}_n(t_f) \) is the final state of a control-free trajectory, while the term \( \vecbf{x}(t_f)\) is the final state of a trajectory under the influence of the control input.
We use the matrix \( Q = \text{diag} \bracket{1~\>0~\> 1~\> 1~\>0~\>1 } \in \R^{6 \times 6}\) to represent the mapping onto \( \Sigma \).
Maximization of the distance between \( \vecbf{x}_n \) and \(\vecbf{x} \), on the \Poincare section defined in~\cref{eq:poincare_section}, is equivalent to the minimization of \( J \) defined in~\cref{eq:cost}.
We ensure that the trajectories intersect the \Poincare section through the use of terminal constraints.
In addition, we use the terminal constraints to define a specific direction along which we seek to minimize the cost~\cref{eq:cost}.
Since the \Poincare section is four-dimensional, we parameterize a direction in \( \R^4 \)  using three angles \( \phi_1, \phi_2 , \phi_3 \).
The terminal constraints are given in terms of these angles as
\begin{align}\label{eq:terminal_constraints}
    \begin{split}
        m_1 &= y = 0 , \\
        m_2 &= \parenth{\sin \phi_{1_{d}}} \parenth{ x_1^2 + x_2^2 + x_3^2 + x_4^2} - x_1^2 = 0, \\
        m_3 &= \parenth{\sin \phi_{2_{d}}} \parenth{ x_2^2 + x_3^2 + x_4^2} - x_2^2 = 0, \\
        m_4 &= \parenth{\sin \phi_{3_{d}}} \parenth{ 2 x_3^2 + 2 x_3 \sqrt{x_4^2 + 2 x_4^2}} - x_3 - \sqrt{x_4^2 + x_3^2} = 0 ,
    \end{split}
\end{align}
where we make use of the difference states \( \parenth{x_1, x_2 ,x_3, x_4 }\) defined as
\begin{align}\label{eq:diff_states}
    \begin{split}
        x_1 &= x(t_f) - x_n(t_f) , \\
        x_2 &= z(t_f) - z_n(t_f) , \\
        x_3 &= \dot{x}(t_f) - \dot{x}_n(t_f) , \\
        x_4 &= \dot{z}(t_f) - \dot{z}_n(t_f) . \\
    \end{split}
\end{align}
We select the terminal time, \( t_f \), from the time required for the uncontrolled trajectory to return back to the \Poincare section.
The constraint \( m_1 = 0 \) ensures that the terminal state lies on the \Poincare section.
The constraints \( m_2, m_3, m_4 \) are used to define a direction on the \Poincare section.
\Cref{eq:diff_states} represents the difference between our controlled and uncontrolled trajectory on the \Poincare section.
We approximate the entire reachable set by discretization  over the space of angles \(\phi_1, \phi_2, \phi_3 \).
By convention we assume that the angles lie in the following range
\begin{align*}
    \phi_1, \phi_2 &\in [ 0, \pi ) ,\\
    \phi_3 &\in [ 0 , 2 \pi ) ,
\end{align*}
such that we parameterize all directions on the three sphere, \(\S^3\).
Finally, we also incorporate the control acceleration magnitude constraint as
\begin{align}\label{eq:control_constraint}
    c(\vecbf{u}) = \vecbf{u}^T \vecbf{u} - u_m^2 \leq 0 ,
\end{align}
where \( u_m \) is the maximum acceleration possible by the propulsion system.
This constraint assumes that the control acceleration may be orientated in any direction yet the acceleration magnitude is variable but bounded.
The goal is to determine the control history \( \vecbf{u}(t) \) such that the cost function~\cref{eq:cost} is minimized while subject to the equations of motion~\cref{eq:body_eoms} and the constraints~\cref{eq:control_constraint,eq:terminal_constraints}.

% How do you solve the optimal control problem
We apply a standard calculus of variations approach to solve our optimal control problem~\cite{bryson1975}.
Using the Euler-Lagrange equations we arrive at the necessary conditions for optimality
\begin{align}\label{eq:necc_conditions}
    \begin{split}
        \dot{\vecbf{x}} ^T &= \deriv{H}{\vecbf{\lambda}} ,\\
        \dot{\vecbf{\lambda}}^T &= \deriv{H}{\vecbf{x}} , \\
        0 &= \deriv{\phi}{x}^T + \deriv{\vecbf{m}}{x}^T \vecbf{\beta} - \vecbf{\lambda}^T(t_f) , \\
        0 &= \deriv{H}{\vecbf{u}} + \mu^T \deriv{c}{\vecbf{u}} ,
    \end{split}
\end{align}
where the Hamiltonian, \( H\), is defined as
\begin{align}\label{eq:hamiltonian}
    H = \vecbf{\lambda}_r^T \vecbf{v} + \vecbf{\lambda}_v^T \parenth{\vecbf{g}(\vecbf{r}) + \vecbf{h}(\vecbf{v}) + \vecbf{u}}.
\end{align}
The costate is given by \( \vecbf{\lambda} = \bracket{ \vecbf{\lambda}_r~\> \vecbf{\lambda}_v }^T \in \R^{6 \times 1}\). 
The vector \( \vecbf{\beta} \in \R^{4 \times 1} \) are the additional Lagrange multiplers associated with the terminal constraints in~\cref{eq:terminal_constraints}, and \( \mu \) is a Lagrange multipler associated with the control constraint in~\cref{eq:control_constraint}.

We can redefine the optimal control in terms of the costate by rewriting the necessary condition as
\begin{align*}
    \vecbf{u} = - \frac{u_m^2}{2 \mu} \vecbf{\lambda}_v .
\end{align*}
We use this along with the control constraint to solve for the Lagrange multiplier \( \mu \)
\begin{align*}
    \mu = \pm \frac{u_m}{2} \norm{\vecbf{\lambda}_v} .
\end{align*}
Finally, we use the second-order necessary condition to determine the correct sign of \( \mu \) and find the optimal control input for the reachable set as
\begin{align}\label{eq:optimal_control}
    \vecbf{u} = - u_m \frac{\vecbf{\lambda_{\vecbf{v}}}}{\norm{\vecbf{\lambda_{\vecbf{v}}}}} .
\end{align}

This optimal control formulation results in a two point boundary value problem.
We use a shooting method to determine the initial costates, \( \vecbf{\lambda}(t_0)\), such that the terminal constraints are satisfied.
In addition, we implement a multiple shooting method which sub-divides the the entire trajectory into small sub-intervals~\cite{stoer2013}.
The multiple shooting method reduces the sensitivity of the terminal states to the initial conditions.
Using a shorter time interval alleviates many of the issues of single shooting approaches, which suffer from convergence difficulties near the optimal solution.
To ensure a continuous trajectory we incorporate additional constraints
\begin{align}\label{sec:interior_constraints}
    \begin{split}
        \vecbf{x}(t_{m-1}^{-}) &= \vecbf{x}(t_{m}^{+}) , \\
        \vecbf{\lambda}(t_{m-1}^{-}) &= \vecbf{\lambda}(t_{m}^{+}),        
    \end{split}
\end{align}
which ensure that both the state and costate are continuous at the patch point between segment \( m-1 \) and \( m\).

\section{Numerical Simulation}\label{sec:simulation}

We present an example transfer of a spacecraft about the asteroid 4769 Castalia. 
Our equations of motion, given by~\cref{eq:eoms}, are an idealized version of the dynamics of a spacecraft.
For example, the model does not include the effect of mass transfer from propellant usage. 
We instead model the control input as a generic acceleration vector in the body-fixed asteroid frame. 
The acceleration magnitude constraint in~\cref{eq:control_constraint} is chosen to emulate a physically realizable thruster system.
In this analysis, we assume \( u_m = \SI{0.1}{\milli\meter\per\second\squared}\) which is equivalent to a thrust of approximately \SI{100}{\milli\newton} for a \SI{1000}{\kilo\gram} spacecraft.
This amount of thrust is typical of many current ion or hall effect thrusters~\cite{goebel2008 ,choueiri2009}.

\begin{figure}[htbp]
    \centering 
    \begin{subfigure}[htbp]{0.45\textwidth} 
        \includegraphics[width=\textwidth]{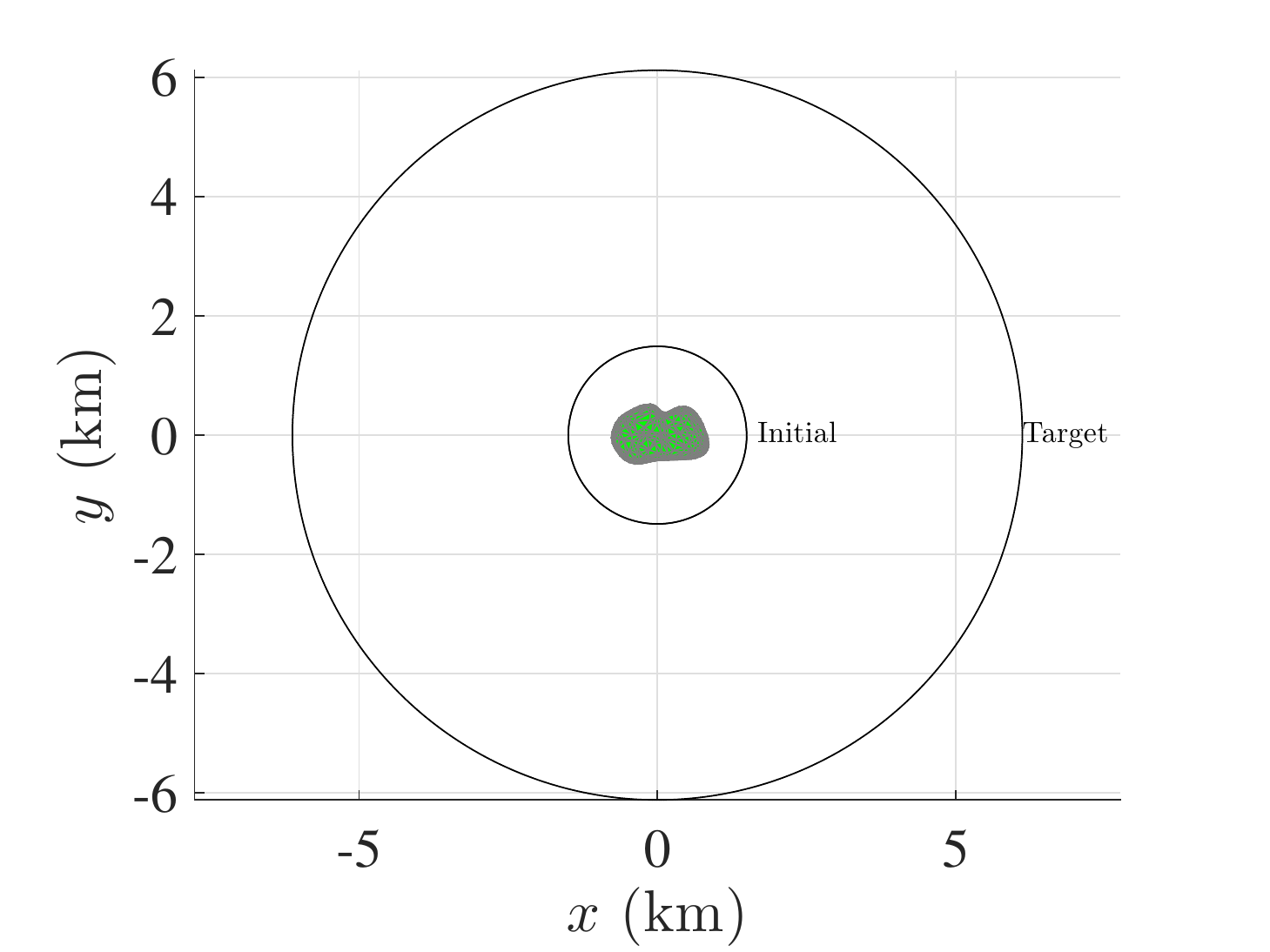} 
        \caption{Equatorial View} \label{fig:eq_initial_transfer} 
    \end{subfigure}~ %add desired spacing between images, e. g. ~, \quad, \qquad, \hfill etc. %(or a blank line to force the subfigure onto a new line) 
    \begin{subfigure}[htbp]{0.45\textwidth} 
        \includegraphics[width=\textwidth]{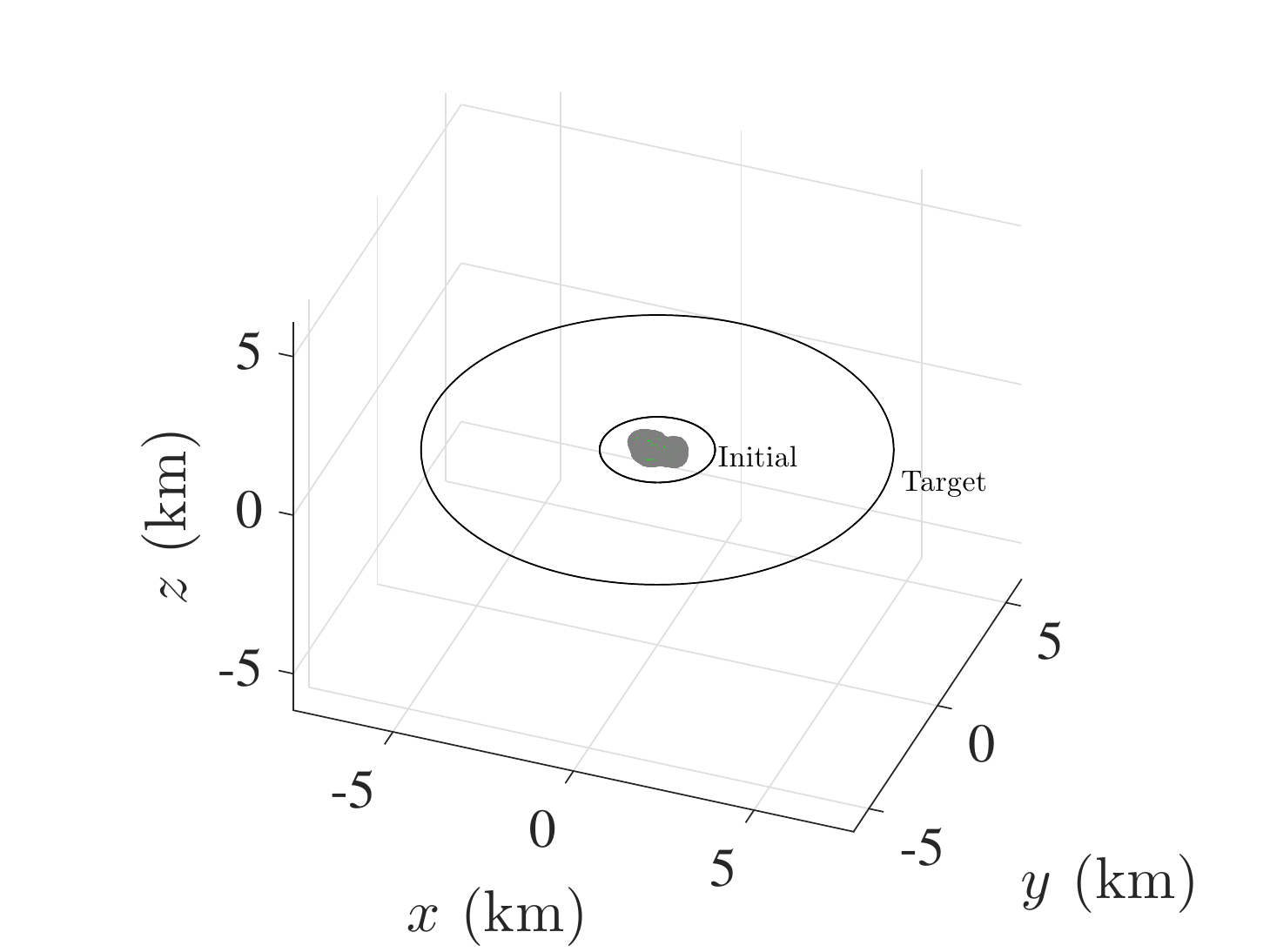} 
        \caption{3D view} \label{fig:initial_transfer_3d} 
    \end{subfigure} ~ %add desired spacing between images, e. g. ~, \quad, \qquad, \hfill etc. %(or a blank line to force the subfigure onto a new line) 
    \caption{Initial and target periodic orbits}
    \label{fig:initial_transfer} 
\end{figure}

The objective is to transfer the spacecraft between two periodic equatorial orbits about Castalia.
The initial and target orbits are periodic solutions about Castlia computed using the method introduced by Reference~\citenum{scheeres2003}.
The initial conditions for both orbits are defined in the body-fixed frame as
\begin{align}\label{sec:initial_transfer}
    \vecbf{x}_i = 
    \begin{bmatrix}
        1.4973 \\ 0 \\ 0.0061 \\ 0\\ -0.0009 \\ 0
    \end{bmatrix} ,
    \quad
    \vecbf{x}_t =
    \begin{bmatrix}
        6.1175 \\ 0 \\ 0.0001 \\ 0\\ -0.0025 \\ 0
    \end{bmatrix} .
\end{align}
\Cref{fig:initial_transfer} shows the initial, \( \vecbf{x}_i \), and target, \( \vecbf{x}_t\), periodic orbits which lie in the equatorial plane of Castalia.
Our goal is to transfer from a lower altitude to a higher altitude while remaining in the equatorial plane of the asteroid.
This type of scenario would occur frequently during mapping and observation missions to asteroids.
In this transfer example we also have used a reduced model of Castalia.
Rather than using the full \num{4092} face model we reduce the number of faces to \num{1024} using the Matlab function \verb+reducepatch+. 
This greatly speeds up the computation with only a small difference in the gravitational field.

We first compute the reachability set originating from the initial periodic orbit at \( \vecbf{x}_i\) for a fixed time of flight and bounded control magnitude as defined previously.
The reachability set is computed by solving the two-point boundary value problem using a multiple shooting algorithm to satisfy the necessary conditions in~\cref{eq:necc_conditions}.
The reachability set is generated on the lower dimensional \Poincare section and is composed of the terminal states in the \( \parenth{x,z,\dot{x},\dot{z} } \) space.
We approximate the reachability set by discretization of each of the angles \( \phi_1, \phi_2 , \phi_3 \) into ten discrete steps. 
This results in a total of \(10^3\) trajectories which approximate the reachability set on the \Poincare section.

\begin{figure}[htbp]
    \centering 
    \begin{subfigure}[htbp]{0.45\textwidth} 
        \includegraphics[width=\textwidth]{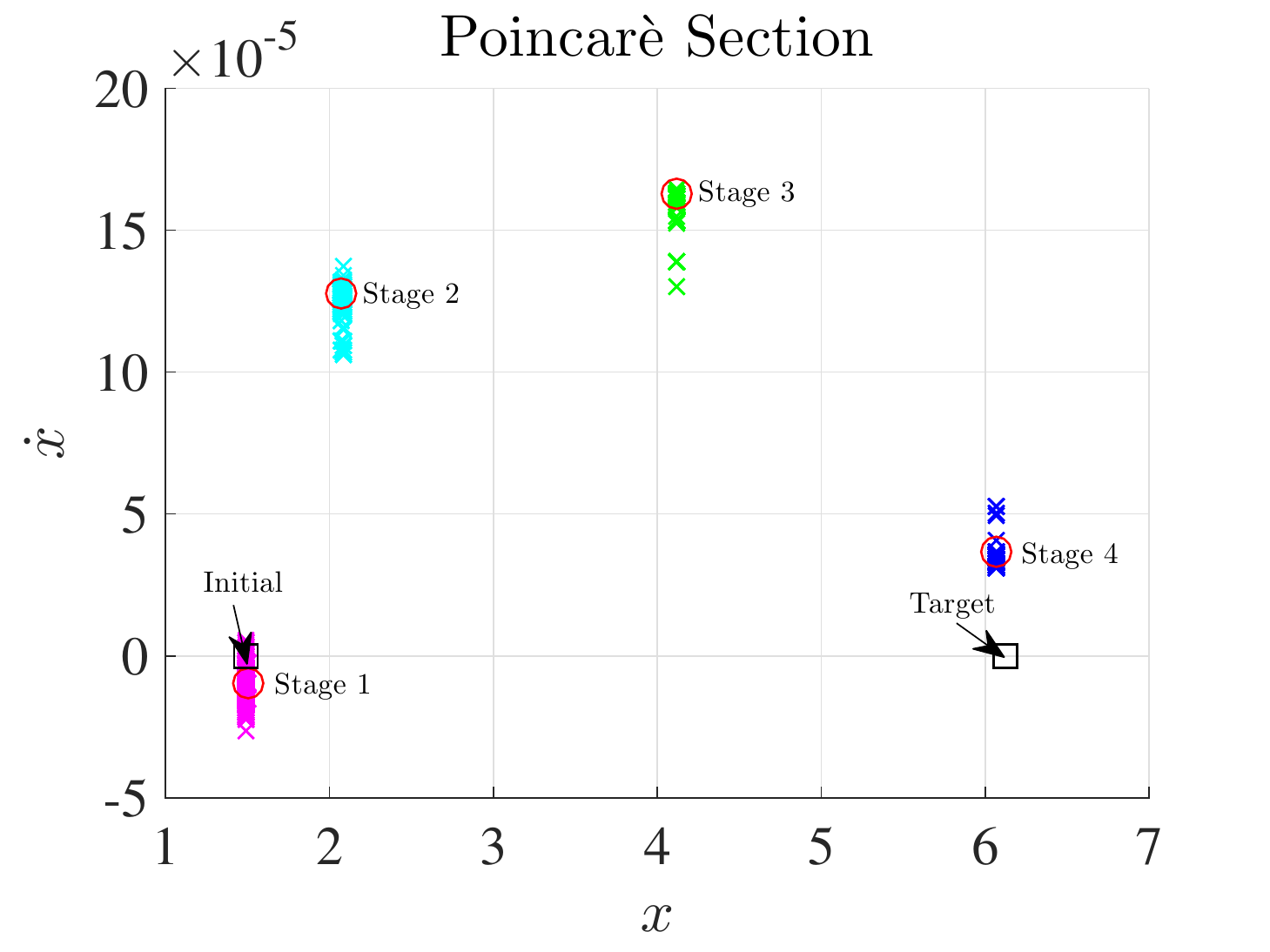} 
        \caption{\( x \) vs. \( \dot{x} \) \Poincare section} \label{fig:poincare_xvsxdot} 
    \end{subfigure}~
    \begin{subfigure}[htbp]{0.45\textwidth} 
        \includegraphics[width=\textwidth]{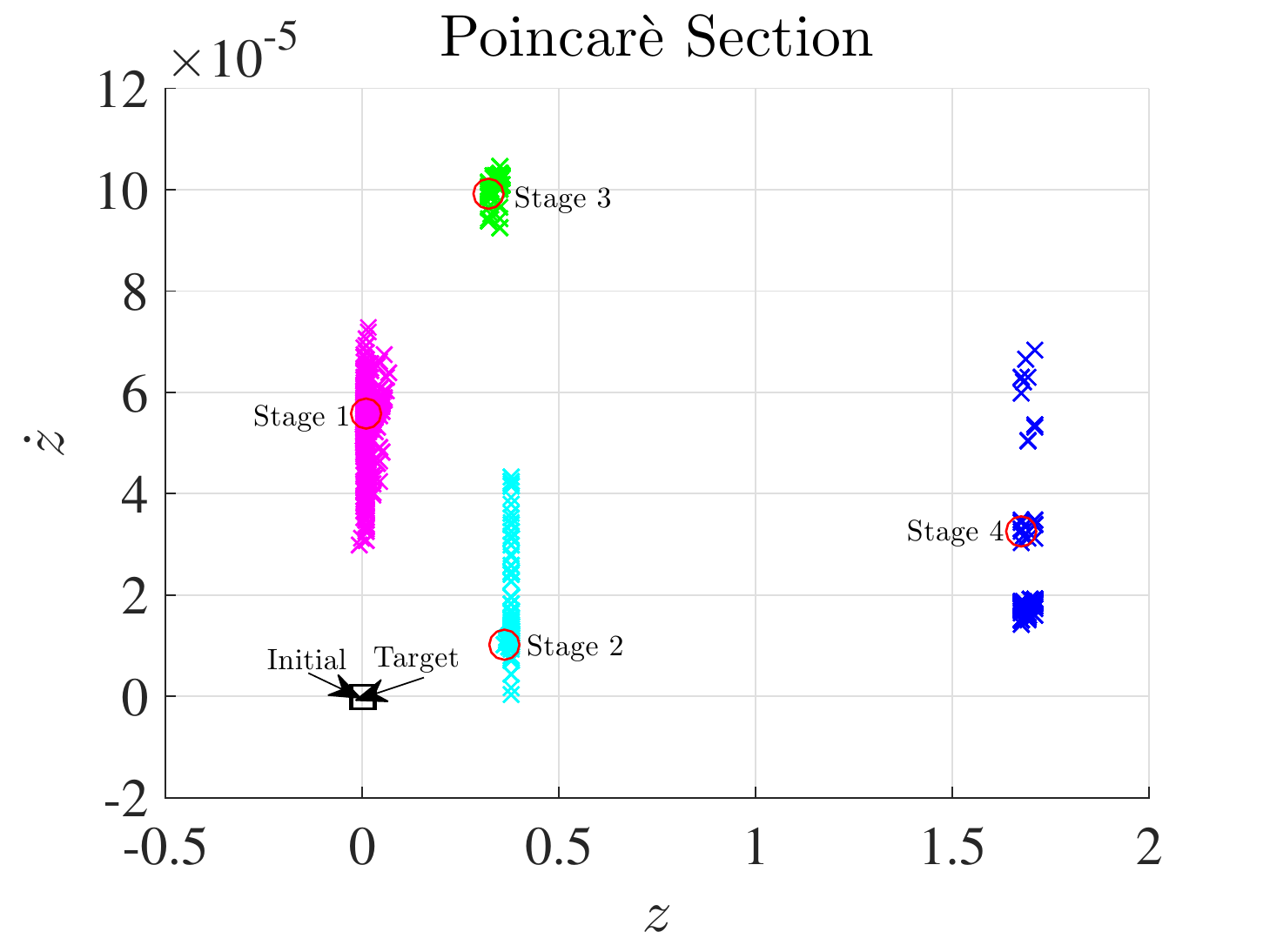} 
        \caption{\( z \) vs. \( \dot{z} \) \Poincare section} \label{fig:poincare_zvszdot} 
    \end{subfigure}
    \caption{\Poincare section visualization \label{fig:poincare_section}}
\end{figure}

\begin{figure}[htbp] 
    \centering 
    \begin{subfigure}[htbp]{0.3\textwidth} 
        \includegraphics[width=\textwidth]{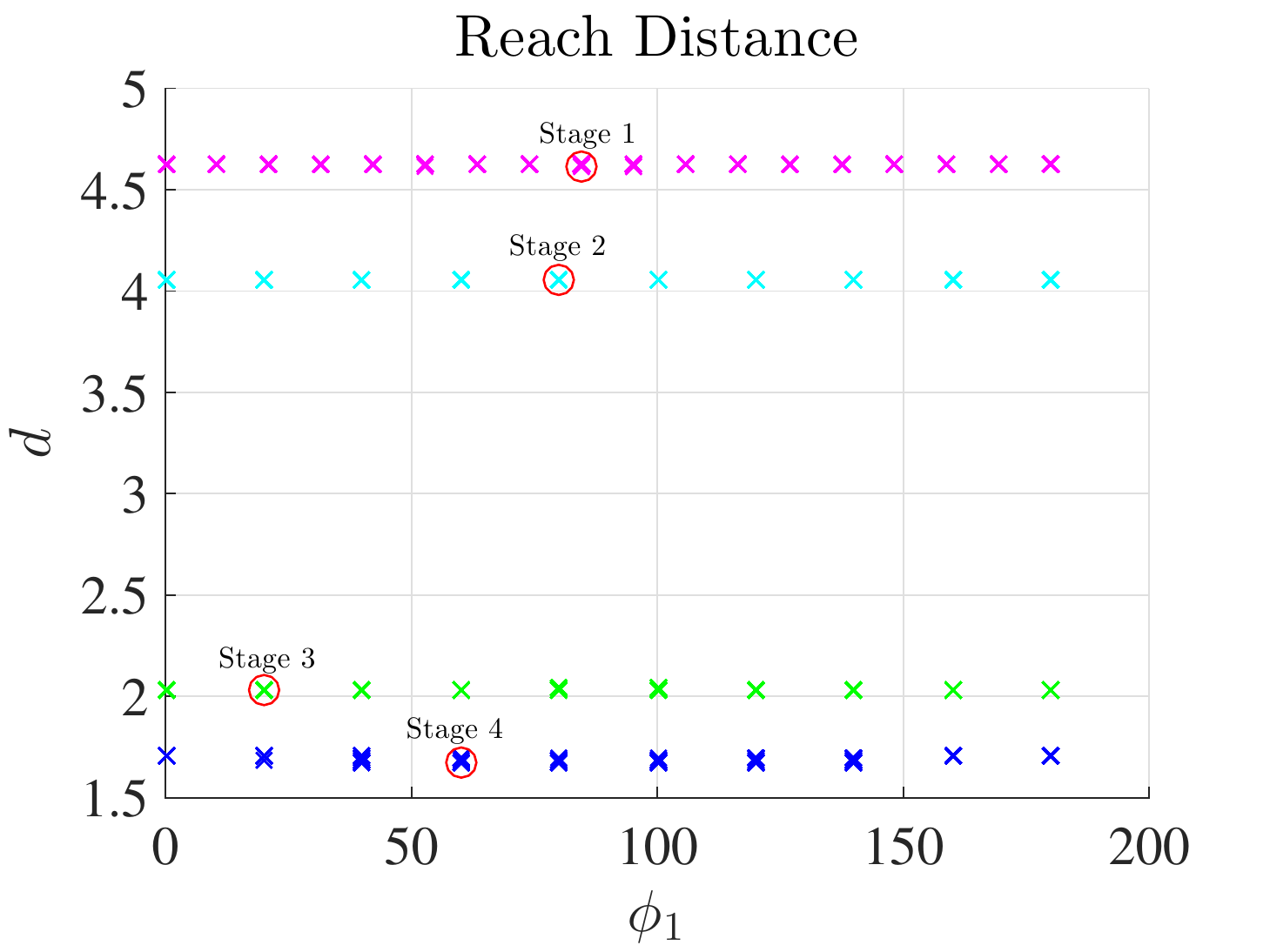} 
        \caption{ \( \phi_1 \)} \label{fig:phi1} 
    \end{subfigure}~
    \begin{subfigure}[htbp]{0.3\textwidth} 
        \includegraphics[width=\textwidth]{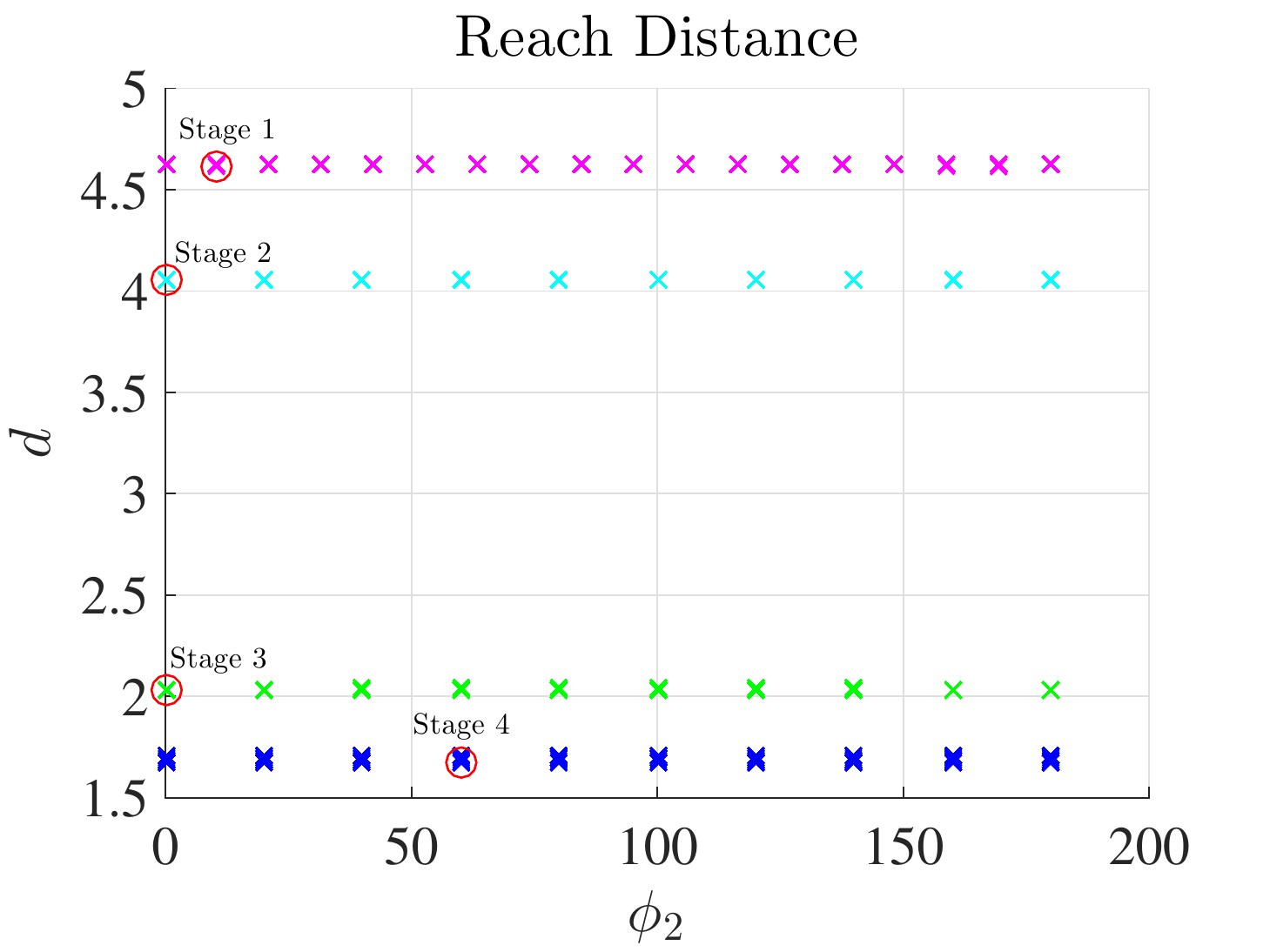} 
        \caption{\( \phi_2 \)} \label{fig:phi2} 
    \end{subfigure}~
    \begin{subfigure}[htbp]{0.3\textwidth} 
        \includegraphics[width=\textwidth]{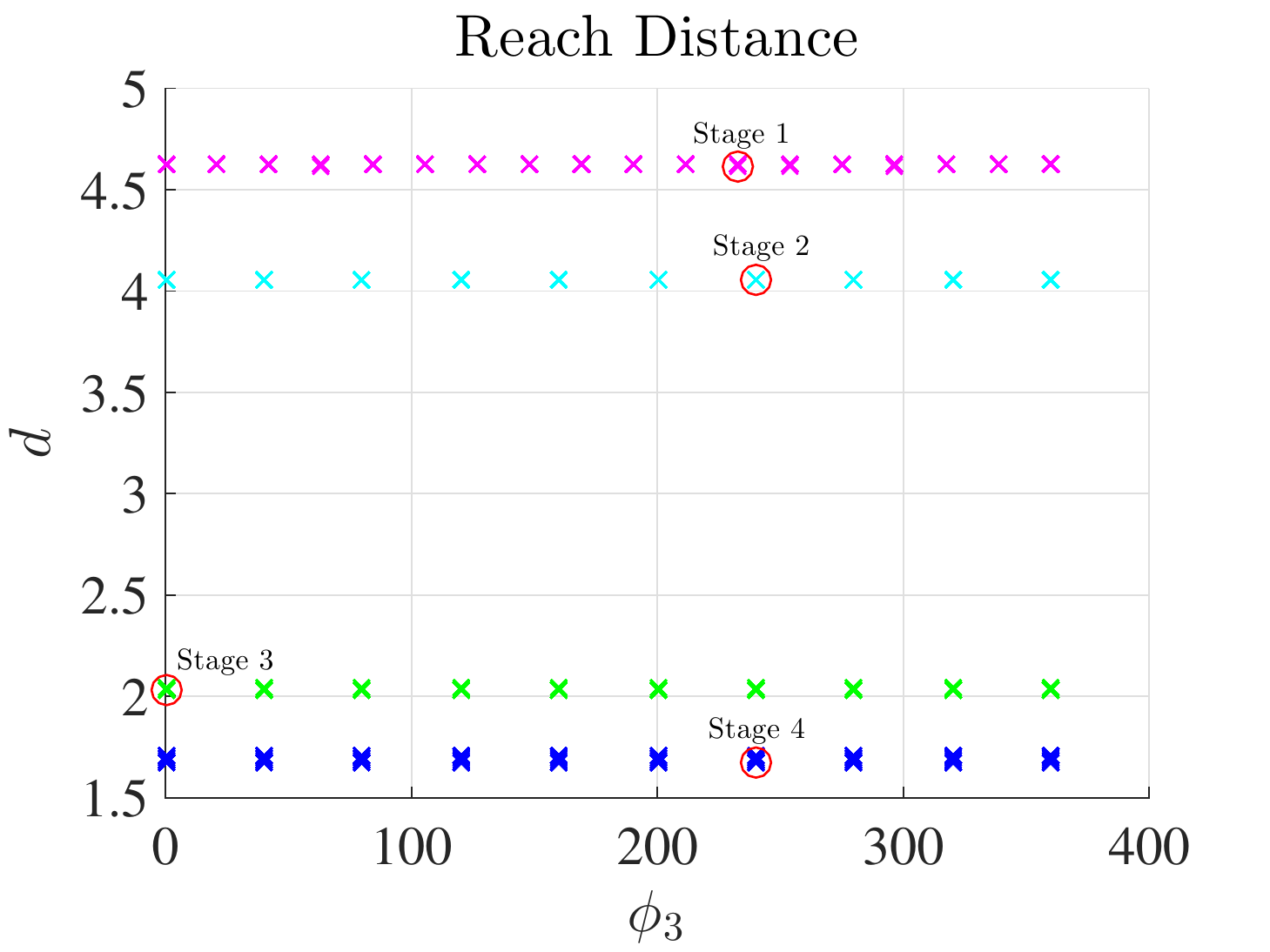} 
        \caption{\( \phi_3 \)} \label{fig:phi3} 
    \end{subfigure} 
    \caption{Variation of \(d(\vecbf{x}_f,\vecbf{x}_t)\) due to \( \phi_i\)}
    \label{fig:phi_distance} 
\end{figure}

We visualize the section using the two figures in~\cref{fig:poincare_section}.
These two-dimensional sections allow us to visualize the four-dimensional \Poincare section defined by~\cref{eq:poincare_section}.
The first stage of the transfer is represented by the magenta markers in~\cref{fig:poincare_section}.
From~\cref{fig:poincare_xvsxdot}, we can see that the reachability set has grown in the \( \dot{x} \) dimension but has not been enlarged much in the \( x \) direction.
Similarly,~\cref{fig:poincare_zvszdot} shows an increase in the \( \dot{z} \) component.
From the reachability set, we chose a trajectory and terminal state which minimizes a distance metric \( d(\vecbf{x}_f,\vecbf{x}_t) \) to the desired target
\begin{align}\label{eq:reach_dist}
    d = \sqrt{k_x \parenth{x_f - x_t }^2 + k_z \parenth{z_f - z_t }^2 + k_{\dot{x}}\parenth{\dot{x}_f - \dot{x}_t }^2 + k_{\dot{z}}\parenth{\dot{z}_f - \dot{z}_t }^2} ,
\end{align}
where \( k_x, k_z, k_{\dot{x}}, k_{\dot{z}} \) are used to weight the relative importance of each of the components of the \Poincare section.
\Cref{fig:phi_distance} shows the distance to the target for the chosen discretization of \( \phi_i \).

The trajectory which minimizes \( d \) is indicated by the red marker in~\cref{fig:poincare_section,fig:phi_distance}.
Since the first reachability set does not include the target we use the minimum state from the first stage to initialize another reachability computation.
Once again we compute the reachability set by discretization of the angles \( \phi_i \) on the \Poincare section.
This second stage, represented by the cyan markers in~\cref{fig:poincare_section,fig:phi_distance}, further increases the \( x, z\) components but does not reach the target orbit.
As a result, a third and forth stage are generated in a similar manner and shown by the green and blue markers in~\cref{fig:poincare_section,fig:phi_distance}, respectively.
We can see in~\cref{fig:poincare_section} that the reachability set of the forth stage includes both the \( x \) and \( z\) components of the target periodic orbit.
At the same time there is a relatively large difference between the \( \dot{x} \), \( \dot{z} \), and \( z \) components of the forth stage and the target orbit.
In practice this is not a large concern as we have direct control over the spacecraft velocity via the control input and the equatorial plane still remains within the reachability set of the transfer.

\begin{figure}[htbp] 
    \centering 
    \begin{subfigure}[htbp]{0.45\textwidth} 
        \includegraphics[width=\textwidth]{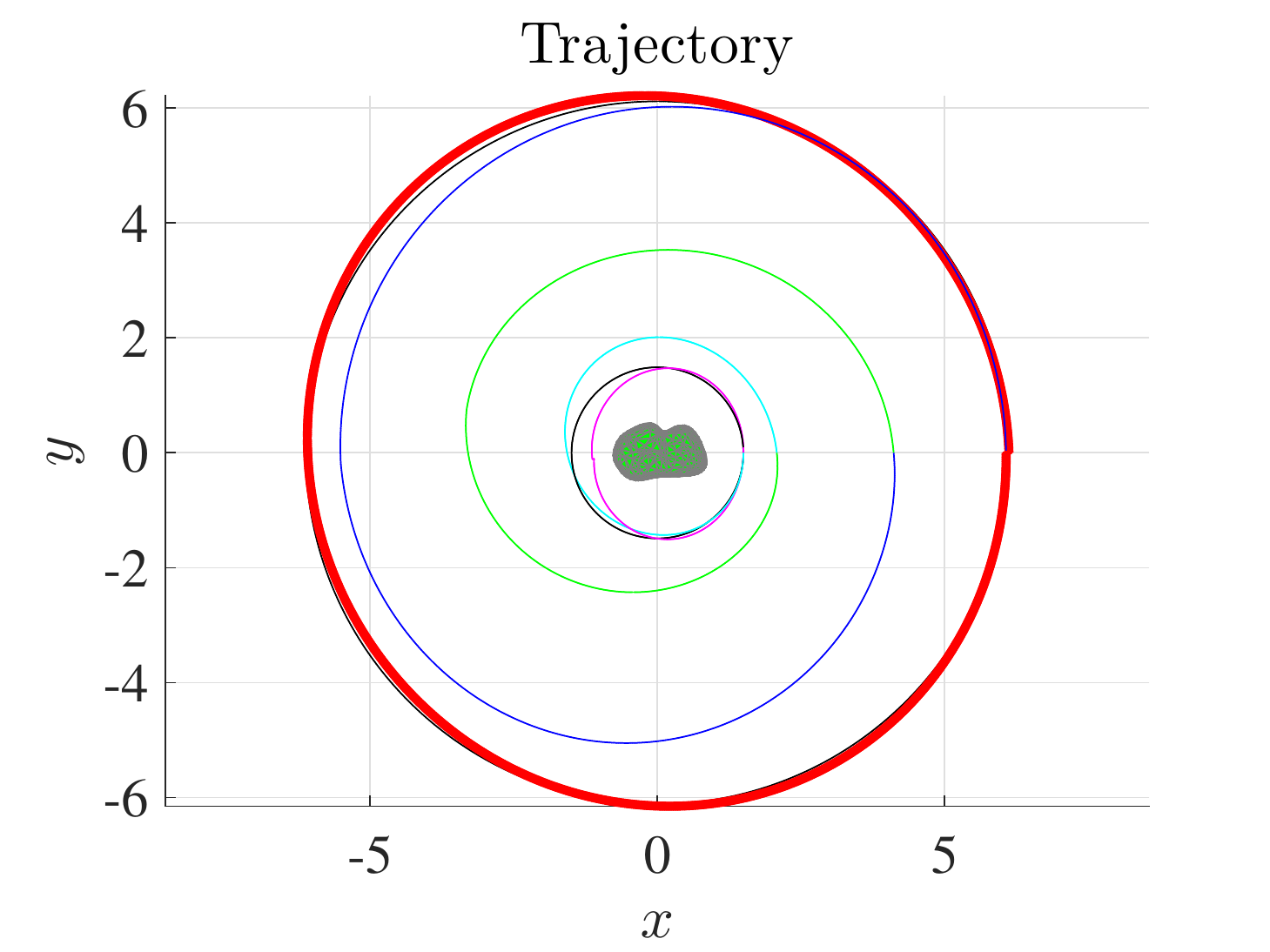} 
        \caption{Equatorial view of transfer} \label{fig:trajectory_up} 
    \end{subfigure}~
    \begin{subfigure}[htbp]{0.45\textwidth} 
        \includegraphics[width=\textwidth]{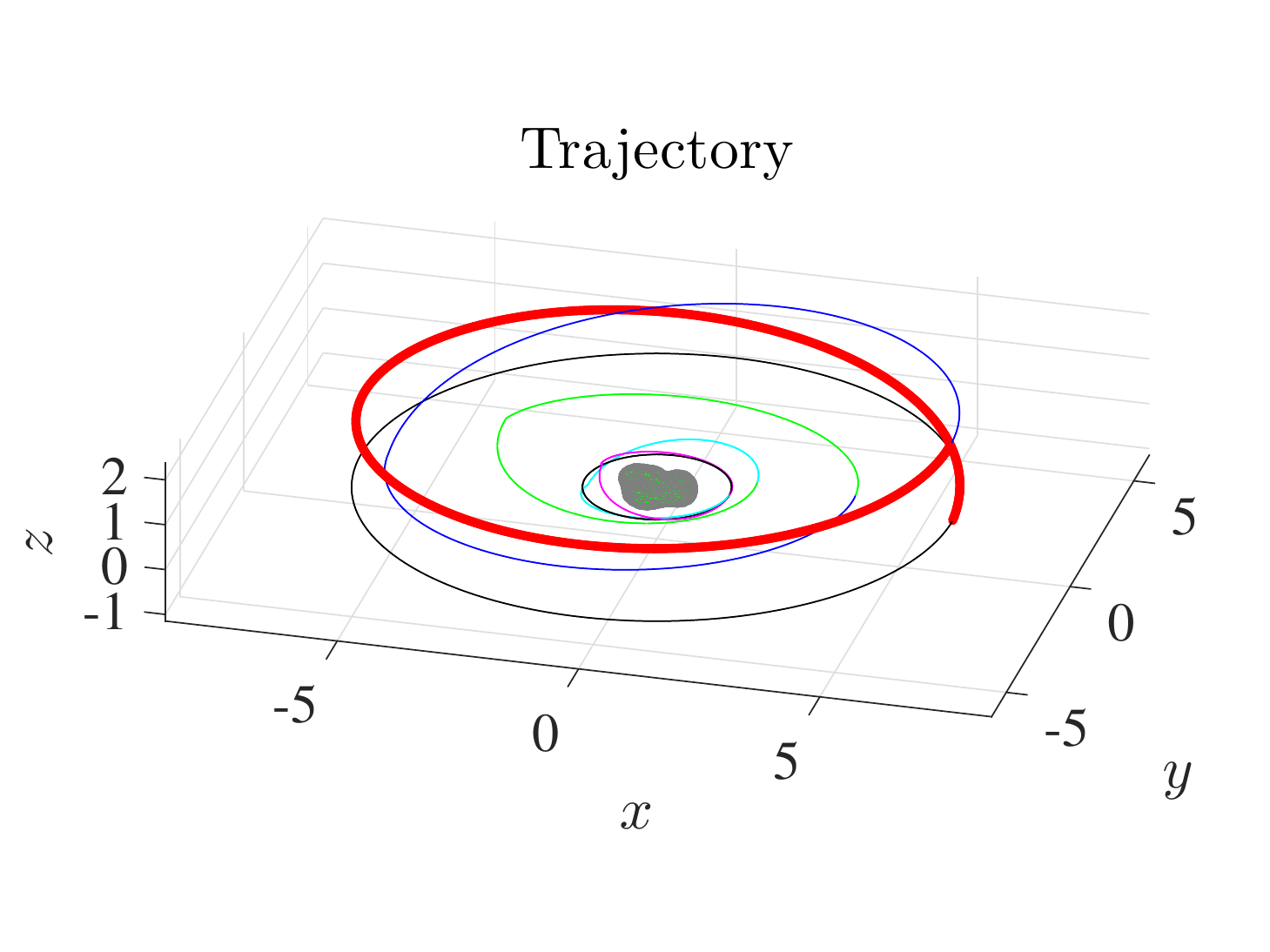} 
        \caption{Out of plane view} \label{fig:trajectory_3d} 
    \end{subfigure}~ 
    \caption{Complete transfer trajectory}
    \label{fig:trajectory} 
\end{figure}

With the reachability set encompassing the target orbit, we can generate a final transfer to the target.
We use the minimum state calculated from the final stage to serve as the initial condition of the transfer.
A final optimal transfer is computed to satisfy the fixed terminal state \( \vecbf{x}(t_f) = \vecbf{x}_t \) and the bounded control magnitude constraint.
\Cref{fig:trajectory} shows the entire transfer trajectory, from the four reachable set trajectories as well as the final transfer to the target.
It is interesting to note that while both the initial and target periodic orbit lie in the equatorial plane, the reachability trajectories show a relatively large out of plane component during the transfers.
In spite of this out of plane movement, the reachability set approaches and meets the target orbit. 
\begin{figure}
    \centering
    \includegraphics[width=0.45\textwidth]{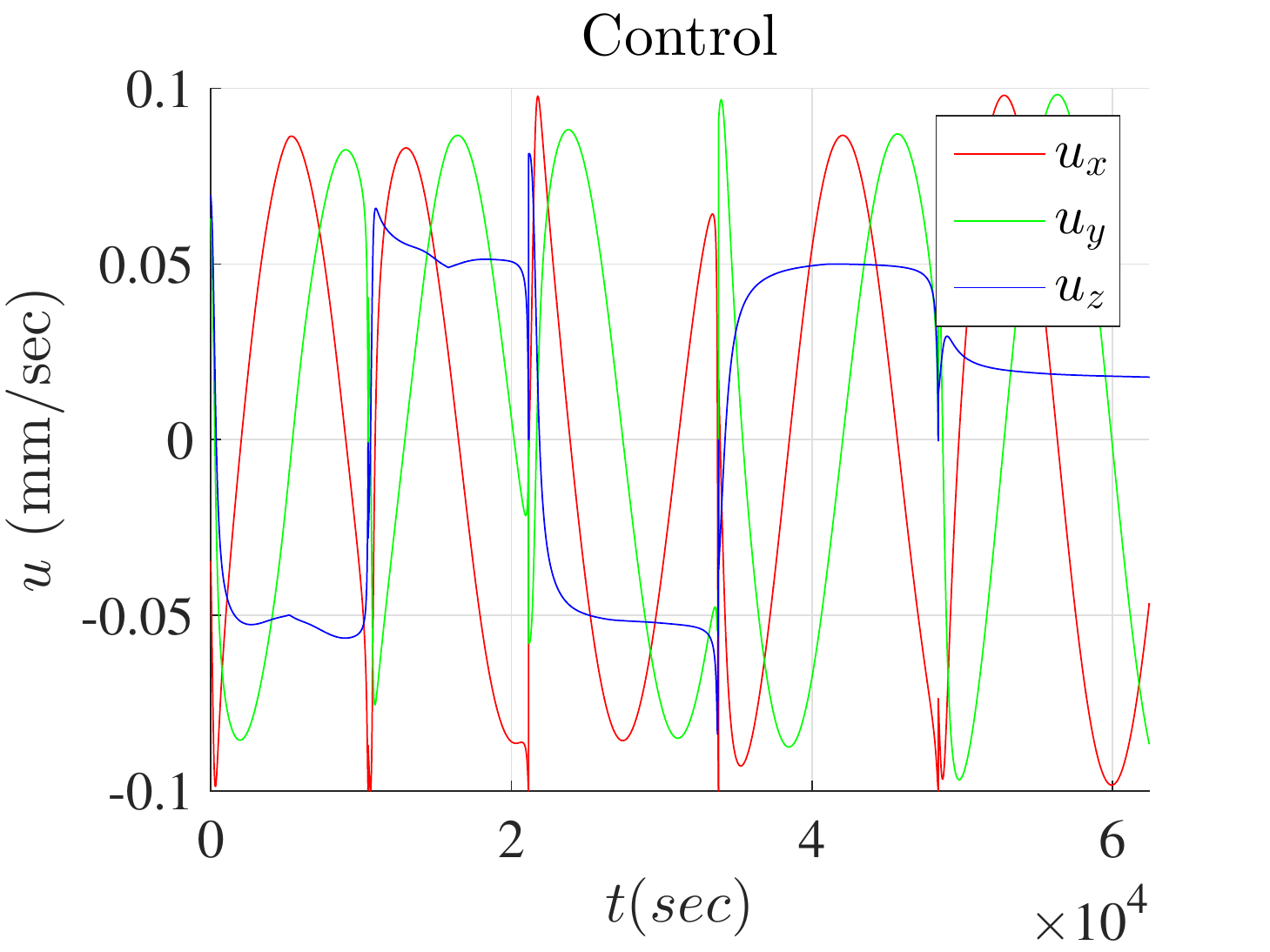}
    \caption{Control history \label{fig:control}}
\end{figure}
\Cref{fig:control} shows the control input required during the maneuver.
We can see that the control constraint in~\cref{eq:control_constraint} is satisfied over the entire trajectory.

\section{Conclusions}\label{sec:conclusions}

In this paper, an optimal transfer process, which combines the concepts of reachability sets on a \Poincare section, is used to generate a transfer between periodic orbits about the asteroid 4769 Castalia.
We have linked several computations of the reachability set on a \Poincare section in order to design a transfer trajectory.
The use of the \Poincare section allows for trajectory design on a lower dimensional space and is an extension of its well-known use in the analysis of periodic orbits.
We use an indirect optimal control formulation to incorporate a control magnitude constraint and several terminal state constraints.
Utilizing the reachability set alleviates the need to determine accurate initial conditions that allow for convergence of the optimal solution.
In addition, the use of the polyhedron gravitational model gives simple method of extending this work to any small body that also possesses a defined shape model.

\section*{Acknowledgments}\label{sec:acknowledgments}

This research has been supported in part by NSF under the grants CMMI-1243000, CMMI-1335008, and CNS-1337722.

\bibliographystyle{aiaa}
\bibliography{library}

\end{document}